\pdfoutput=1
\documentclass{article}
\usepackage[margin=25mm]{geometry}
\usepackage{amsmath}
\usepackage{amsfonts}
\usepackage{amssymb}
\usepackage{graphicx}

\usepackage[utf8]{inputenc} % правильне кодування
\usepackage[T2A]{fontenc}
\usepackage[english]{babel}

\usepackage{xcolor}
\usepackage[all]{xy}

\newtheorem{theorem}{Theorem}

\title{Topological structure of Morse functions on the projective plane}
\author{Svitlana Bilun, Alexandr Prishlyak, Serhii Stas and Alina Vlasenko}

\date{\today}

\begin{document}
\maketitle

\begin{abstract}
To investigate the topological structure of Morse functions on the projective plane we use the Reeb graphs. 
We describe it properties and prove that it is a complete topological invariant of simple Morse function on $\mathbb{R} P^2$. We prove recurent formulas for the number of Reeb graphs with the given number of sadles (vertex).
\end{abstract} \hspace{10pt}



\section*{Introduction}

In the papers of Kronrod \cite{Kronrod1950} and Reeb \cite{Reeb1946}, an invariant (the Reeb graph) was introduced, which describes the topological properties of functions on two-dimensional manifolds. In the case of simple Morse functions on closed oriented two-dimensional manifolds, it is a complete topological invariant under fiberwise equivalence. If a linear order is specified on the set of vertices, then it becomes a complete topological invariant under topological equivalence. In the case of a non-orientable two-dimensional manifold \cite{lychak2009morse}, as well as manifolds with boundary \cite{hladysh2017topology, hladysh2019simple}, additional information is needed to construct a complete topological invariant . For non-compact manifolds \cite{prishlyak2002morse}, the Reeb graph can be a non-Hausdorff space, which causes additional difficulties when working with it. If the function is not simple, then in addition to the Reeb graph, information about the structure of the function in the neighborhood of each critical level is needed \cite{Bolsinov2004}.

Another way to specify the topological structure of Morse functions is to use Morse-Smale vector fields with given function values at singular points \cite{lychak2009morse, Smale1961}. Therefore, the topological classification of Morse--Smale vector fields is closely related to the classification of functions.

Topological properties of smooth function on closed surfaces was also investigated in \cite{prishlyak2001conjugacy, hladysh2019simple, hladysh2017topology,  prishlyak2002morse, Prishlyak2000,  prishlyak2007classification, lychak2009morse, prish2002Morse, prish2015top, prish1998sopr, Bilun2013def, Bilun2002,  Sharko1993}, on surfaces with the boundary in \cite{Hladysh2016, hladysh2019simple} and on closed 3-manifolds in  \cite{prishlyak1999equivalence, prishlyak2003regular}. 
%Bilun13,

We recomend \cite{Kybalko2018, Oshemkov1998, Peixoto1973, prishlyak1997graphs, prishlyak2020three, Akchurin2022, Prishlyak2022, prislyak2017morse,  kkp2013,  prish2002vek,     Prishlyak2021,  Prishlyak2020, Kybalko2018,  Prishlyak2019} on classifications of flows on closed 2- manifolds and 
\cite{Losieva2017, prislyak2017morse, Prishlyak2022, prishlyak2003sum, Prishlyak2003, Prishlyak2019} on manifolds with the boundary.
Topological properties of flows on 3-manifolds was investigated in \cite{prish1998vek,  prish2001top, Prishlyak2002beh2, prishlyak2002morse1,  prishlyak2003regular, Prishlyak2002b, prishlyak2005complete, Prishlyak2007, Hatamian2020, BPP2022, Bilun2022}.

The main invariants of graphs and their embeddings in surfaces can be founded in \cite{prishlyak1997graphs, Harary69, pontr86, tatt88, HW68, GT87}.

The topological properties of the projective plane are described in many geometry textbooks, but we would also recommend \cite{Bilun22Projective}.

The main purpose of this paper is to study the topological properties of the Reeb graphs of simple Morse functions on the projective plane and to count their number for a fixed number of critical points.

\section{Morse function on $RP^2$}

Let M be a smooth closed two-dimensional manifold, $f : M \to R$ be a smooth function. 

A critical point $x \in M$ is called non-degenerate, if the matrix $H = \left(
\frac{\partial^2 f(x)}{\partial x_i \partial x_j} 
\right)$ 
in some local coordinates $x_1, x_2$ is non-degenerate. On a two-dimensional manifold there are 3 types of non-degenerate critical points: a local minimum, a saddle and a local maximum. 

A smooth function $f : M \to R$ is called \textit{a Morse function}, if all its critical points are non-degenerate. A Morse function is called \textit{simple} if all its critical points lie on different levels, $f(p) \ne f(q)$, if $p \ne q$. 

A component of the level line $f^{ -1} (y)$ of the Morse function is called \textit{a fiber}. 

Two Morse functions are called \textit{fiber equivalent} if there is a homeomorphism of the surface onto itself which maps fibers of one function to fibers of another one, and the local minima to the local minima. %A neighborhood of the critical fiber which is foliated into level lines of the function and considered to within the fiber equivalence is called an atom. We consider only simple Morse functions. 

\begin{figure}[ht]
\center{ \includegraphics[height=4.5cm]{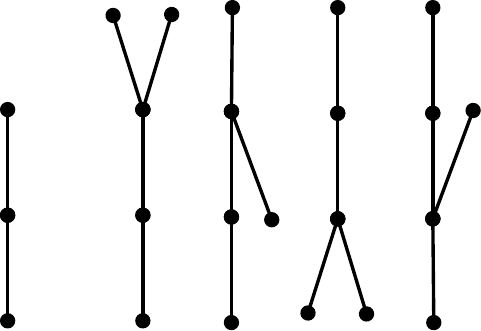} }
\caption{1 Reeb graphs   with one saddle and 4  Reeb graphs   with two saddles}
\label{Morse12}
\end{figure}

\begin{figure}[ht]
\center{ \includegraphics[height=14.0cm]{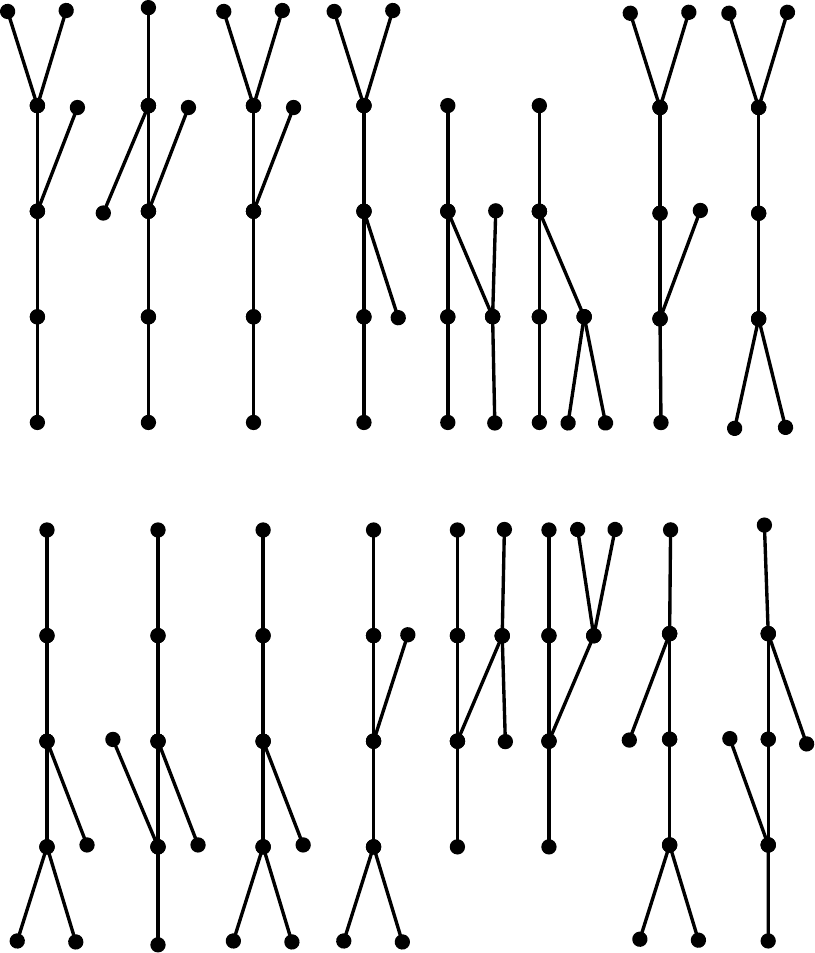} }
\caption{16 Reeb graphs of Morse functions on $RP^2$  with three saddle critical points}
\label{Morse3}
\end{figure}

The quotient space $M/ \sim$ with orientation of edges according to the direction of the increase of the function is called a Reeb graph, where $f : M \to R$ is a Morse function, $x_1 \sim x_2$ if $x_1$ and $x_2$ belongs to one fiber. Reeb graphs are considered to within the isomorphism of oriented graphs. %The atom can be one of three types, — trousers, inverted trousers, and a nonorientable atom. The first atom corresponds to a vertex of the Reeb graph with two edges one of which is directed toward and the other outwards the vertex. The second atom corresponds to a vertex of the Reeb graph with one edge directed towards and two edges directed outwards the vertex. The nonorientable atom corresponds to the vertex of valency 2. 
Examples of Reeb graphs are shown on Fig. \ref{Morse12}, \ref{Morse3}.

Two simple Morse functions on the orientable surface are fiber equivalent if and only if their Reeb graphs are isomorphic ([1], Theorem 2.4, p. 71). On the non-orientable surfaces there are firer non-equivalent functions with the same Reeb graphs. So additional invariants are nessersity to clasify Morse functions. 
But in the case of projective plane we have the following

\begin{theorem}
The Reeb graph of a simple Morse function on the projective plane has the following properties:

1) the graph is a tree;

2) the graph has one vertex of degree 2, the other vertices have degree 1 or 3;

3) only vertices of degree 1 are sinks and sources.
\end{theorem}
\textbf{Proof.} 1) If the graph is not a tree, then it has a cycle, which corresponds to a connected sum with a torus or a Klein bottle, which means that the genus of a non-orientable surface is greater than or equal to 2, which is impossible.

2) Each vertex of degree 2 corresponds to a non-orientable atom (a Mjöbius strip with a hole). If there is more than one such vertex, then the genus of the surface will be greater than one.

3) sinks and sources correspond to local extrema of the function, and vertices of degree 2 and 3 correspond to saddle points.

\begin{theorem}
Two simple Morse functions on $\mathbb{R}P^2$ are fiber equivalent if and only if their Reeb graphs are isomorphic. 
\end{theorem}
\begin{figure}[ht]
\center{ \includegraphics[height=5.0cm]{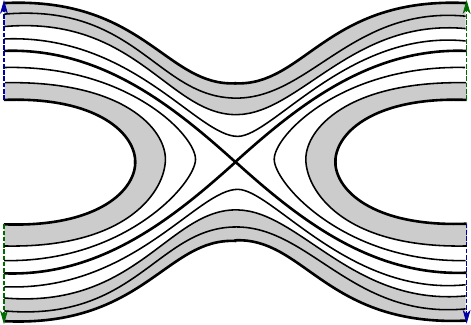} }
\caption{An neighborhood of the critical level}
\label{atom}
\end{figure}
\textbf{Proof}.
Consider a neighborhood $U$ of the critical level corresponding to a vertex of valence 2 on the Reeb graph. After cutting along two regular trajectories of the gradient field, it will take the form shown in Fig. \ref{atom}. In addition to this neighborhood, the function is a Morse function on an oriented surface, so the isomorphism of the Reeb graphs on them defines the topological equivalence of the restrictions of the functions to these parts. Let us continue these homeomorphisms on $U$. To do this, we divide Y into two parts: one part is the collar (neighborhood of the edge), and the second is the central one (complement to the collar). If the constructed topological equivalences preserve orientation on each of the components of the boundary Y, then on the central part we define an identical mapping, if they reverse on each of them, then central symmetry, if they reverse orientation on one of the edges, then axial symmetry about the horizontal or vertical axis. Since any homeomorphism of a circle is isotopent to either the identity mapping or symmetry about the axis, then these isotopies define topological equivalence of functions on the collar. By construction, all three mappings (on the collar, the central part, and complements to $U$) coincide at the edges, then they define the required topological equivalence.

\begin{figure}[ht]
%\center{ \includegraphics[height=6.0cm]{Morse8-1RP2.pdf} }
%\caption{Reeb graphs of Morse functions on $RP^2$  with 4 saddle critical points. Part 1}
%\label{Morse8-1}
%\end{figure}
%\begin{figure}[ht]
\center{ \includegraphics[height=10.0cm]{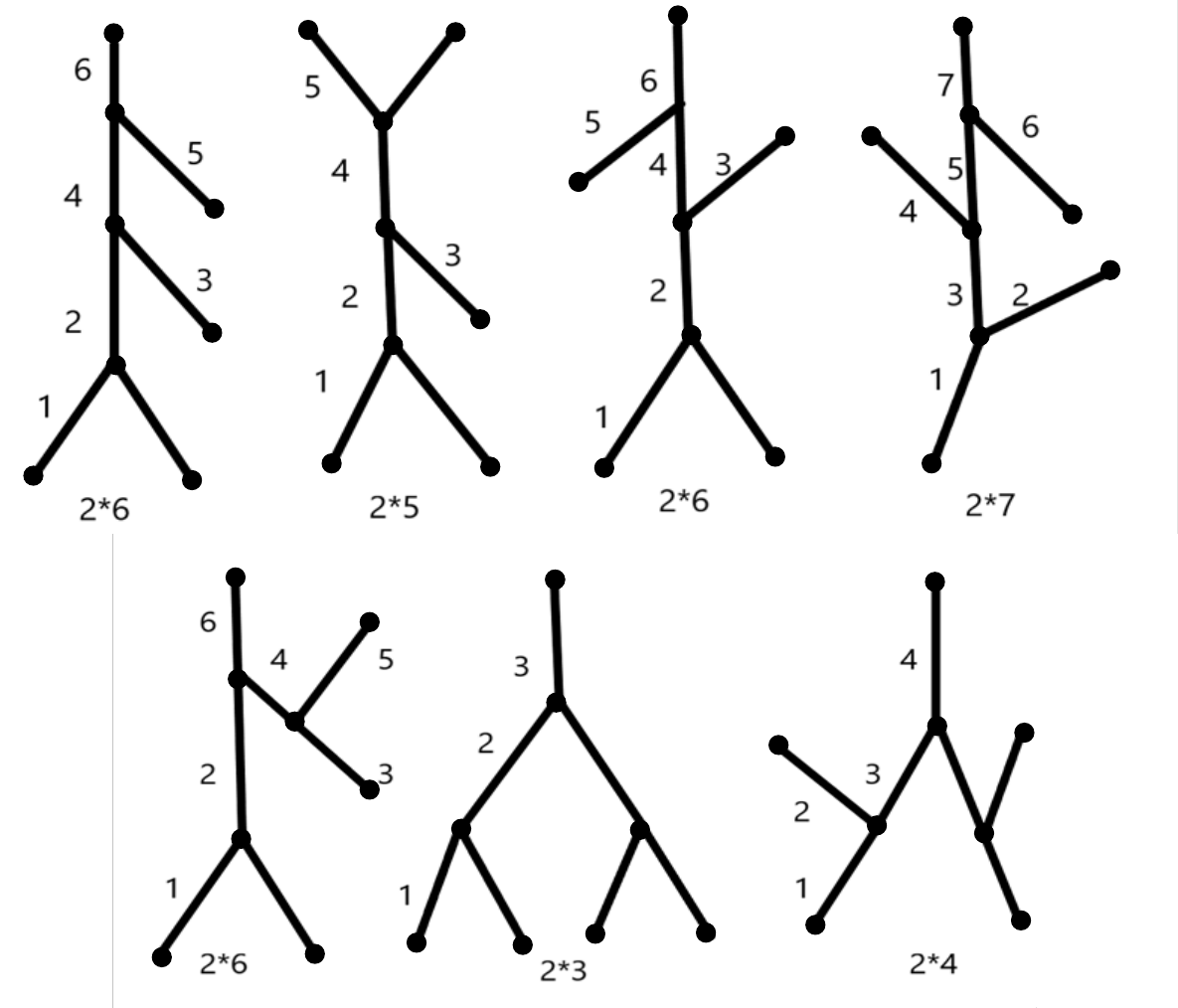} }
\caption{74 Reeb graphs of Morse functions on $RP^2$  with four saddle critical points%. Part 2
}
\label{Morse4}
\end{figure}
\section{Number of Reeb graphs}
To obtain the number of topologically nonequivalent Morse functions on the projective plane, we first calculate the number of rooted Reeb graphs with a given number of saddles.
\begin{theorem}
The number $K_k$ of root oriented Reeb graphs with $k$ saddles can be calculated using the following recursive formulas: for even $k=2n$
$$K_{2n}=3(K_0 K_{2n-1}+K_1 K_{2n-2}+ \ldots +K_{n-1}K_n);$$
for odd $k=2n+1$
$$K_{2n+1}=3(K_0 K_{2n}+K_1 K_{2n-1}+ \ldots +K_{n-1}K_{n+1})+\frac{3K_n^2+K_n}{2}.$$
\end{theorem}
\textbf{Proof.}
Consider the saddle closest to the root. It is a vertex of degree 3. Remove this saddle from the graph. We obtain three connected components. One of them is an edge connecting the saddle with the root. The other two components can be thought of as rooted Reeb graphs, where the root is the remote saddle. Let the first of these graphs have at least vertices than the second. Then three options are possible: 

1) in both graphs, the edge incident to the root is directed upwards in the original graph, 

2) in the first graph this edge is directed upwards, and in the second downwards, 

3) in the second graph it is directed upwards, and in the first graph it is directed down. 

In cases where the specified edge is directed downward, it causes the orientations on the rooted tree to be reversed. Note that both graphs are interconnected only by the condition on the total number of vertices, which must be equal to the number of vertices in the original graph. If the number of vertices in both graphs is different, then they cannot be isomorphic. Therefore, the total number of non-isomorphic original graphs with such subgraphs will be equal to the product of the numbers of non-isomorphic subgraphs. So we get the first formula. If the number of vertices in the first graph is equal to the number of vertices in the second graph, then the second and third options coincide and must be counted once. In addition, in the first variant, not isomorphic graphs were counted twice in the product, but isomorphic ones once. Considering all this, we obtain the second formula.

Using these formulas, we obtain such initial values.

$K_0=1$

$K_1= 2$

$K_2= 6$

$K_3= 25$

$K_4= 111$

$K_5= 540$

$K_6= 2736$

$K_7= 14396$

$K_8= 77649$

$K_9= 427608$

$K_{10}= 2392866$

$K_{11}= 13570386$

$K_{12}= 77815161$

$K_{13}= 450418536$

$K_{14}= 2628225684$

\begin{theorem}
The number $N_k$ of topologically non-equivalent simple Morse functions on $\mathbb{R}P^2$ can be calculated by the following formula 
$$N_k=K_0 K_{k-1} + K_1 K_{k-2}+ \ldots +K_{k-1}K_{0}.$$
\end{theorem}
\textbf{Proof.} The Reeb graph of a simple Morse function on the projective plane is partitioned by a vertex of degree 2 into two rooted oriented Reeb graphs. Therefore, the total number of Reeb graphs is equal to the sum of the products of the corresponding rooted Reeb graphs. 

Using it, we obtain

$N_1=1$

$N_2= 4$

$N_3= 16$

$N_4= 74$

$N_5 = 358$

$N_6 = 1824$

$N_7 = 9589$

$N_8 = 51766$

$N_9 = 285035$

$N_{10} = 2178244$

$N_{11}= 9046744$

$N_{12} = 51876774$

$N_{13} = 300278112$

$N_{14} = 1752150456$

$N_{15} = 10295599780$

All possible Reeb graphs with one and two saddles are shown in Fig.\ref{Morse12}, with three saddles in Fig.\ref{Morse3}, and with four edges in Fig.\ref{Morse4}.

\section*{Conclusion} 

We proved that Reeb graph is complete topological invariant of simple Morse function on $\mathbb{R} P^2$. 
A recurrent formula was derived for counting the number of topologically non-equivalent functions. The results of calculations for functions with no more than 15 saddles (31 critical points) are given.
We hope that the obtained results can be generalized to other surfaces.

%\bibliographystyle{plain}
%\bibliography{prishd}

\textsc{Taras Shevchenko National University of Kyiv}

Svitlana Bilun \ \ \ \ \ \ \  \textit{Email address:} \text{ bilun@knu.ua}   \ \ \ \ \ \ \ \ \ 
\textit{ Orcid ID:}  \text{0000-0003-2925-5392}

Alexandr Prishlyak \ \textit{Email address:} \text{ prishlyak@knu.ua} \ \ \ \
\textit{ Orcid ID:} \text{0000-0002-7164-807X}

Serhii  Stas \ \ \ \ \   \ \textit{Email address:} \text{ stasserhiy380@gmail.com} \ 
\textit{ Orcid ID:} \text{0009-0006-1241-7497}

Alina  Vlasenko  \ \ \textit{Email address:} \text{ vlasenko\_aline@ukr.net} \ 
\textit{ Orcid ID:} \text{0009-0000-0213-959X}

\end{document}